\documentclass{article}
\usepackage[utf8]{inputenc}

\RequirePackage{amsthm,amsmath,amsfonts,amssymb}

\usepackage{amsmath}
\RequirePackage[numbers]{natbib}
\RequirePackage[ruled,vlined,linesnumbered]{algorithm2e}
\RequirePackage[colorlinks,citecolor=blue,urlcolor=blue]{hyperref}
\RequirePackage{graphicx,comment}

\newtheorem{theorem}{Theorem}
\newtheorem{proposition}[theorem]{Proposition}
\newtheorem{lemma}[theorem]{Lemma}

\newtheorem{definition}[theorem]{Definition}
\newtheorem{problem}[theorem]{Problem}

\newtheorem{example}[theorem]{Example}
\newtheorem{conjecture}[theorem]{Conjecture}

\newcommand{\tconv}{\text{tconv\,}}
\newcommand{\Trop}{{\text{trop}}}
\DeclareMathOperator*{\argmin}{arg\,min}

\newcommand{\RR}{\mathbb{R}}
\newcommand{\R}{\mathbb{R}}

\usepackage{tabularx, url}

\newcommand{\Proj}{\text{Proj\,}}

\title{Tropical Data Science}
\author{Ruriko Yoshida}
\date{}

\begin{document}

\maketitle
\section*{Outline}
\begin{enumerate}
   \item Introduction
    \item Data Science Overview
    \begin{enumerate}
        \item Supervised Learning (Inferential Statistics) 
        \begin{enumerate}
            \item Classifications
            \item Regression 
        \end{enumerate}
        \item Unsupervised Learning (Descriptive Statistics)
    \end{enumerate}
    \item Phylogenetics to Phylogenomics 
    \begin{enumerate}
        \item Phylogenetic Trees
        \item Space of Phylogenetic Trees
    \end{enumerate}
    \item Basics in Tropical Geometry
    \item Tropical Unsupervised Learning (Tropical Descriptive Statistics)
    \begin{enumerate}
        \item Tropical Fermat Weber Points
        \item Tropical Fr\'ecet Means
        \item Tropical Principal Component Analysis (PCA)
    \end{enumerate}
    \item Tropical Supervised Learning (Tropical Inferential Statistics)
    \begin{enumerate}
        \item Tropical Classifications
        \begin{enumerate}
            \item Tropical Support Vector Machines
            \item Tropical Linear Discriminant Analysis
        \end{enumerate}
        \item Tropical Regression
    \end{enumerate}
\end{enumerate}

\begin{abstract}
    Phylogenomics is a new field which applies to tools in phylogenetics to genome data.  Due to a new technology and increasing amount of data, we face new challenges to analyze them over a space of phylogenetic trees.  Because a space of phylogenetic trees with a fixed set of labels on leaves is not Euclidean, we cannot simply apply tools in data science.  In this paper we survey some new developments of machine learning models using tropical geometry to analyze a set of phylogenetic trees over a tree space. 
\end{abstract}

{\em Key Words:} Machine Learning Models, Max-Plus Algebra, Phylogenetics, Phylogenomics, Tropical Geometry, Ultrametrics.

\section{Introduction}

Due to increasing amount of data today, data science is one of most exciting fields.  It finds applications in statistics, computer science, business, biology, data security, physics, and so on. 
Most of statistical models in data sciences assume that data points in an input sample are distributed over a Euclidean space if they have numerical measurements.  However, in some cases this assumption can be failed.  For example, a {\em space of phylogenetic trees} with a fixed set of leaves is an union of lower dimensional cones over $\mathbb{R}^e$, where $e = {N \choose 2}$ with $N$ is the number of leaves \cite{AK}.  Since the space of phylogenetic trees is an union of lower dimensional cones, we cannot just apply statistical models in data science to a set of phylogenetic trees \cite{YZZ}.  

There has been much work in spaces of phylogenetic trees.  In 2001, Billera-Holmes-Vogtman (BHV) developed a notion of a space of phylogenetic trees with a fixed set of labels for leaves \cite{BHV}, which is a set of all possible phylogenetic trees with the fixed set of lebels on leaves and is an union of orthants, each orthant is for all possible phylogenetic trees with a fixed tree topology.  In this space, two orthants are next to each other if the tree topology for one orthant is one nearest neighbor interchange (NNI) distance to the tree topology for the other orthant.  They also showed that this space is CAT(0) space so that there is a unique 
shortest connecting paths, or geodesics, between any two points in the space defined by the ${\rm   CAT}(0)$-metric. There is some work in development on machine learning models with the BHV metric.  For example, Nye defined a notion of the first order
principal component geodesic as the unique geodesic with the BHV metric over the the tree space which minimizes the sum of residuals between the geodesic and each data point \cite{Nye}.  However, we cannot use a convex hull under the BHV metric for higher principal components because Lin et. al showed that the convex hull of three points with the BHV metric over the tree space can have arbitrarily high dimension \cite{LSTY}.

In 2004, Speyer and Sturmfels showed a space of phylogenetic trees with a given set of labels on their leaves is a tropical Grassmanian \cite{SS}, which is a tropicalization of a linear space defined by a set of linear equations \cite{YZZ} with the max-plus algebra.  The tropical metric with max-plus algebra on the tree space is known to behave very well \cite{AGNS,CGQ}.  For example, contrarily to the BHV metric, the dimension of the convex hull of $s$
  tropical points is at most $s-1$.

Thus, this paper focuses on the tropical metric over tree spaces.  In this paper we review some development on statistical learning models with the tropical metric with max-plus algebra on tree spaces as well as the tropical projective space, and we overview some open problems.

\section{Data Science Overview}

In this section, we briefly overview statistical models in data science.  For more details, 
we recommend to read {\em Introduction of Statistical Learning with R} \url{http://faculty.marshall.usc.edu/gareth-james/ISL/}.

\begin{figure}
\begin{center}
    \includegraphics[width=\textwidth]{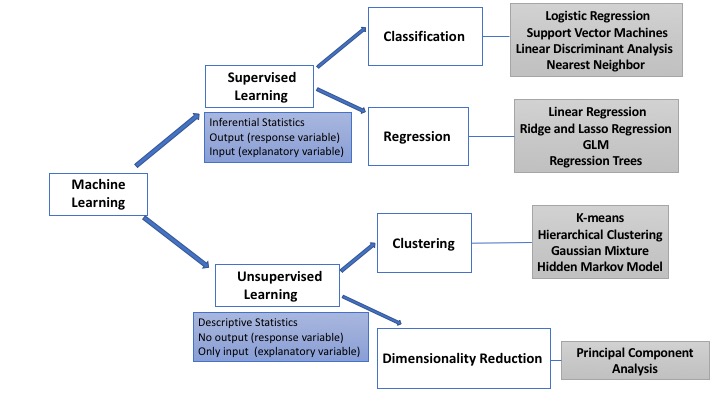}
\end{center}
    \caption{Overview of data science}
    \label{fig:ds}
\end{figure}

In data science there are roughly two sub-branches of data science: unsupervised learning and supervised learning (Figure \ref{fig:ds}).   In unsupervised learning, our goal is to compute a descriptive statistics to see how data points are distributed over the sample space or how data points are clustered together. In statistics, unsupervised learning corresponds to descriptive statistics.   In supervised learning, our goal is to predict/infer the response variable from explanatory variables.  In statistics,  supervised learning corresponds to inferential statistics.  Like unsupervised learning and supervised learning, there are some notations with different names between machine learning and statistics.  Thus we summarize some of the differences in Table \ref{tab:notation}.

\begin{table}[h]
    \centering
    \begin{tabular}{c|c}
         Statistics& Data Science \\ \hline
          descriptive statistics & Unsupervised learning\\
          inferential statistics & Supervised learning \\
         response variable & target variable\\
         explanatory variable & predictor variable\\
         & feature\\
    \end{tabular}
    \caption{There are several notations with different names in statistics and data science.}
    \label{tab:notation}
\end{table}
  
\subsection{Basic Definitions}
\begin{enumerate}
    \item {\em Response variable} -- the variable for an interest in a study or experiment.  It can be called as a dependent variable.  In machine learning it is also called a target variable.
    \item {\em Explanatory variable} -- the variable explains the changes in the response variable.  It can be also called a {\em feature} or independent variable.  In machine learning it is also called feature or predictor.
\end{enumerate}

\subsection{Unsupervised Learning}
Since unsupervised learning is descriptive, there is no response variables. In unsupervised learning, we try to learn how data points are distributed and how they related to each other.  Among them, there are mainly two categories: clustering and dimensionality reduction.  
\begin{itemize}
    \item {\em Clustering} -- grouping data points into subsets by their ``similarity".  These similarities are defined by a user.  These groups are called {\em clusters}. 
    \item {\em Dimensionality reduction} -- reducing the dimension of data points with minimizing the loss of information.  One of the most commonly used methods is principal component analysis (PCA), a dimension reduction procedure via linear algebra.
\end{itemize}

\subsection{Supervised Learning}
Supervised learning is inferential.  Thus, there are the response variable and explanatory variables in an input data set.  Depending on the scale of the response variable, we can separate two groups in supervised learning: classification and regression.  In classification, the response variable has categorical scale and in regression, the response variable has  numerical (interval) scale.  
\begin{itemize}
    \item {\em Classifications} -- the response variable is categorical.  Under classification, there are algorithms like logistic regression, support vector machine, linear discriminant analysis, classification trees, random forests, adaboost and etc.
    \item {\em Regression} -- the response variable is numerical.  There are algorithms like linear regression, regression trees, lasso, ridge regression, random forests,  adaboost and etc.
\end{itemize}

For more details, see the following papers:
\begin{itemize}
    \item  Gareth James, Daniela Witten, Trevor Hastie and Robert Tibshirani. {\rm Introduction of Statistical Learning with R} \url{http://faculty.marshall.usc.edu/gareth-james/ISL/}.
\end{itemize}

\section{Phylogenetics to Phylogenomics}

In this section we overview basics in phylogenetics and basic problem for phylogenomics.  

\subsection{Phylogenetic Trees}

Evolutionary, or phylogenetic trees, show an organism's evolutionary relationships over time, through the use of tree diagrams.
Phylogenetic trees still consist of vertices (nodes) and edges (branches).
Each node in a phylogenetic tree represents a past or present taxon or population: exterior nodes in a phylogenetic tree represent taxon or population at present; and interior nodes represent their ancestors. Edges in a phylogenetic tree have weights and a weight in each edge represents mutation rates multiplied by evolutional time from its ancestor to a taxon.  

\begin{figure}[ht!]
\centering
\includegraphics[width=10cm]{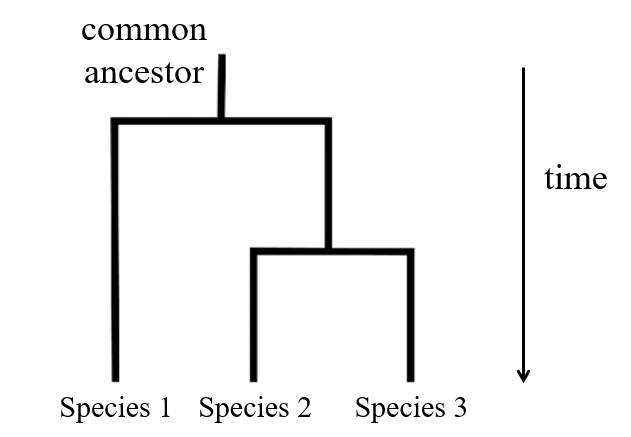}
\caption{Example Rooted Binary Phylogenetic Tree}
\label{fig6}
\end{figure}

In Figure \ref{fig6}, an exterior vertex (leaf or tip) represents the current taxa {\it ($V=\{Species 1, Species 2,Species 3\})$}.
An interior vertex represents an extinct taxa $(V=\{4,5\})$ where ancestors split into two subgroups.
Vertices and edges can still be labeled; however, only the leaves or tips are labeled in a phylogenetic tree.
This is due to the past taxa often being inferred and not exactly known.
Vertices in phylogenetic trees can be DNA sequences, shared genes or interrelated species, depending on the context of the tree.
The root {\it $(V=\{5\})$} of the tree now represents the common ancestor of all leaves, {\it $Species 1$}, {\it $Species 2$}, and {\it $Species 3$}.

Phylogenetic trees are trees.  Thurs, they remain acylic and connected.  In terms of evolutional biology, these properties are intuitive, as a species must evolve from something and also as time progresses species can only evolve forward.
Weights on edges in a tree represent the notion of time.
The distance {\it $(d_{Species1,Species2})$} measures the dissimilarity between {\it Species 1} and {\it Species 2} with respect to time.

Let $N$ be the number of leaves on a phylogenetic tree.    If a total of weights of all edges in a path from the root to each leaf $i \in [N] := \{1, \ldots , N\}$ in a rooted phylogenetic tree $T$ is the same for all leaves $i \in [N]$, then we call a phylogenetic tree $T$ {\em equidistant tree}.  The {\em height} of an equidistant tree is the total weight of all edges in a path from the root to each leaf in the tree. 

\subsubsection{Phylogenetic Tree Reconstruction}

Phylogenetic reconstruction uses genetic data to create an inferential evolutionary (phylogenetic) tree.   
These {\it changing characters} are the mutations in DNA sequences.
DNA sequences represent a shared gene across multiple species.
Trees are excellent at representing the evolutionary changes of this shared gene through node splits and leaves.

Even though we do not discuss details on a phylogenetic tree reconstruction in this paper, multiple steps and techniques are involved in the reconstruction process and 
there are several types of tree reconstruction methods;
\begin{itemize}
    \item {\em Maximum likelihood estimation (MLE) methods} -- These methods describe evolution in terms of a discrete-state continuous-time Markov process.
    \item {\em Maximum Parsimony} -- Reconstructs tree with the
    least evolutionary changes which explain data.
    \item {\em Bayesian inference for trees} -- Use Bayes Theorem and MCMC to estimate the posterior distribution rather than obtaining the point estimation.  
    \item {\em Distance based methods} -- Reconstructing a tree from a distance matrix.
\end{itemize}

\subsection{Space of Phylogenetic Trees}

There are several ways to define a space of phylogenetic trees with different metrics. One of the very well-known tree spaces is Billera-Holmes-Vogtmann tree space. In 2001, Billera-Holmes-Vogtmann (BHV)  introduced a continuous space which models the set of rooted phylogenetic trees with edge lengths on a fixed set
of leaves. In this space, edge lengths in a tree are continuous and we assign a coordinate for each {\em interior edge}.  Note that unrooted trees can be accommodated by designating a fixed leaf node as the root. The {\em BHV tree space} is not Euclidean, but it is non-positively curved, and thus has the property that any two points
are connected by a unique shortest path through the space, called a {\em geodesic}. The distance between two trees is defined as the
length of the geodesic connecting them.  While in this paper, we do not consider the BHV tree space, read \cite{BHV} for interested readers.  

Through this paper, we assume that all phylogenetic trees are {\em equidistant trees}. An equidistance tree is a rooted phylogenetic tree such that the sum of all branch lengths in the unique path from the root to each leaf in the tree, called the {\em height} of the tree, is fixed and they are the same for all leaves in the tree.   In phylogenetics this assumption is fairly mild since the multispecies coalescent model assumes that all gene trees have the same height.
  
  \begin{example}
  Suppose $N = 4$.   Consider two rooted phylogenetic trees with the set of labels on the leaves $\{a, b, c, d\}$ in Figure \ref{fig:equidistant}. Note that for each tree, the sum of branch lengths in the unique path from the root to each leaf is $1$. Therefore  they are equidistant trees with their height are equal to $1$.
    \begin{figure}
      \centering
      \includegraphics[width=0.5\textwidth]{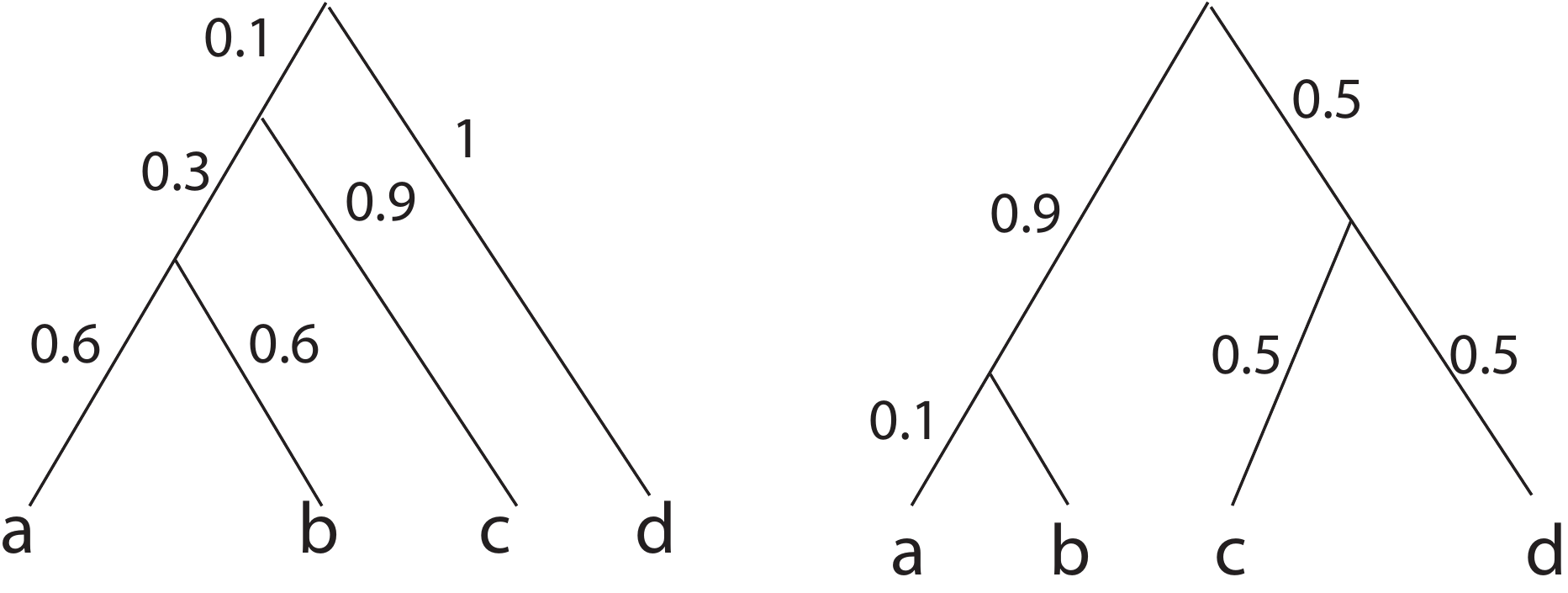}
      \caption{Examples of equidistant trees with $N = 4$ leaves with the set of labels $\{1, b, c, d\}$ and with their height equal to $1$.}
      \label{fig:equidistant}
  \end{figure}
  \end{example}
  
  For the space of equidistant trees with the fixed set of labels on their leaves, the BHV tree space might not be appropriate \cite{treespace:ultra}.  Therefore, we consider the {\em space of ultrametrics}.  To define {\em ultrametrics} and theire relations to equidistant trees, we need to define {\em dissimilarity maps}.  
  
  \begin{definition}\label{def:dissimilarity}[Dissimilarity Map]
  A {\em dissimilarity map} $w$ is a function $w:[N]\times [N]\to \mathbb{R}_{\geq 0}$ such that 
  $$
    w(i,j)=\begin{cases}
    0 & \mbox{if } i = j\\
    \geq 0 &\mbox{else,}
    \end{cases}
  $$ 
    for all $i, \, j\in [N]$. If a dissimilarity map $w$ additionally satisfies the triangle inequality, that is:
\[
w(i,j)\leq w(i,k)+w(k,j),
\]
for all $i, \, j, \, k\in [N]$, then $w$ is called a {\em metric}.
If there exists a phylogenetic tree $T$ such that $w(i, j)$ coincides with the total branch length of the edges in the unique path from a leaf $i$ to a leaf $j$ for all leaves $i, j \in [N]$, then we say $w$ a {\em tree metric}.  If a metric $w$ is a tree metric and $w(i, j)$ is the total branch length of all edges in the path from a leaf $i$ to a leaf $j$ for all leaves $i, j \in [N]$ in a phylogenetic tree $T$, then we say $w$ {\em realises a phylogenetic tree $T$} or $w$ is a {\em realisable} of a phylogenetic tree $T$.
  \end{definition}
 Since
 \[
 w(i, j) = \begin{cases}
 0 & \mbox{if } i = j\\
 w(j, i) & \mbox{else,}
 \end{cases}
 \]
 to simplify we write 
  \[
    w = \left(w(1, 2), w(1, 3), \ldots , w(N-1, N)\right).
    \]
  
   \begin{example}
  We consider equidistant trees in Figure \ref{fig:equidistant}. The dissimilarity map obtained from the left tree in Figure \ref{fig:equidistant} is 
  \[
  (1.2, 1.8, 2, 1.8, 2, 2).
  \]
  Similarly,  the dissimilarity map obtained from the right tree in Figure \ref{fig:equidistant} is 
  \[
  (0.2, 2, 2, 2, 2, 1).
  \]
  Since these dissimilarity maps are obtained from phylogenetic trees, they are tree metrics.
  \end{example}
  
  \begin{definition}[Three Point Condition]
  If a metric $w$ satisfies the following condition: For every distinct leaves $i, j, k \in [N]$,
  \[
    \max\{w(i, j), \, w(i, k), \, w(j, k)\}
    \]
  achieves twice, then we say that $w$ satisfies the {\em three point condition}.
  \end{definition}
  
  \begin{definition}[Ultrametrics]\label{def:ultra}
  If a metric $w$ satisfies the three point condition then $w$ is called an {\em ultrametric}.
  \end{definition}
   
  \begin{theorem}[\cite{JARDINE1967173}]\label{equid:ultra}
  A dissimilarity map $w: [N] \times [N]$ is ultrametric if and only if $w$ is realisable of an equidistant tree with labels $[N]$.  In addition, for each equidistant tree there exists a unique ultrametric.  Conversely, for each ultrametric, there exists a unique equidistant tree. 
  \end{theorem}
  \begin{example}
  We again consider equidistant trees in Figure \ref{fig:equidistant}. The dissimilarity map obtained from the left tree in Figure \ref{fig:equidistant} is 
  \[
  (1.2, 1.8, 2, 1.8, 2, 2).
  \]
  Similarly,  the dissimilarity map obtained from the right tree in Figure \ref{fig:equidistant} is 
  \[
  (0.2, 2, 2, 2, 2, 1).
  \]
  Since these phylogenetic trees are equidistant trees, these dissimilarity maps are ultrametrics by Theorem \ref{equid:ultra}.  
  \end{example}

  From Theorem \ref{equid:ultra} we consider the {\em space of ultrametrics} with labels $[N]$ as a space of all equidistant trees with the label set $[N]$.  Let $\mathcal{U}_N$ be the {space of ultrametrics} for equidistant trees with the leaf labels $[N]$.  
  In fact we can write $\mathcal{U}_N$ as the tropicalization of the linear space generated by linear equations.

  Let $L_N \subseteq \mathbb{R}^e$ be the linear subspace 
  defined by the linear equations such that
  \begin{equation}
  \label{eq:trop_eq}
  x_{ij} - x_{ik} + x_{jk}=0
  \end{equation} 
  for $1\leq i < j < k \leq N$.  
  For the linear equations (\ref{eq:trop_eq}) spanning the linear space $L_N$, the max-plus tropicalization $Trop(L_N)$ of the linear space $L_N$ is the tropical linear space with $w \in \mathbb{R}^e$ such that $$\max\left\{w_{ij}, \, w_{ik}, \, w_{jk}\right\}$$ achieves at least twice for all $i, j, k \in [N]$.  Note that this is exactly the three point condition defined in Definition \ref{def:ultra}.

\begin{theorem}{\cite[Theorem 2.18]{YZZ}}
  \label{tropicalLine}
  The image of $\mathcal{U}_N$ in the tropical projective torus $\mathbb{R}^{n}/\mathbb{R}{\bf 1}$ coincides with $\Trop(L_N)$.  
  \end{theorem}
  
For example, if $N = 4$,  The space of ultrametrics $\mathcal{U}_4$ 
  is a two-dimensional fan with  $15$ maximal cones.  

 For more details, see the following papers:
 \begin{itemize}
     \item  C. Semple and M. Steel. {\rm Phylogenetics}, \cite{Steel2003}.
     \item Lin et al. {\rm Convexity in Tree Spaces} \cite{doi:10.1137/16M1079841}.
 \end{itemize}
 
 \section{Basics in Tropical Geometry}

 Here we review some basics of tropical arithmetic and geometry, as well as setting up the notation through this paper. 

\begin{definition}[Tropical arithmetic operations]
  Throughout this paper we perform arithmetic over the
  max-plus tropical semiring $(\,\mathbb{R} \cup \{-\infty\},\boxplus,\odot)\,$.
  Over this tropical semiring,  the basic tropical
  arithmetic operations of addition and multiplication are defined as the following:
    $$a \boxplus b := \max\{a, b\}, ~~~~ a \odot b := a + b ~~~~\mbox{  where } a, b \in \mathbb{R}\cup\{-\infty\}.$$
Over this tropical semiring, $-\infty$ is the identity element under addition and $0$ is the identity element under multiplication. 
  \end{definition}
  
\begin{example}
 Suppose we have $a = 1, \, b = -3$.  Then 
 \[
 \begin{array}{rcl}
 1 \boxplus (-3) &=& \max\{1, -3\} = 1\\
 1 \odot (-3) &=& 1 + -3 = -2.\\
 \end{array}
 \]
\end{example}
  
  \begin{definition}[Tropical scalar multiplication and vector addition]
  For any $a,b \in \mathbb{R}\cup \{-\infty\}$ and for any $v = (v_1, \ldots ,v_e), w= (w_1, \ldots , w_e) \in (\mathbb{R}\cup-\{\infty\})^e$, tropical scalar multiplication and tropical vector addition are defined as:
    $$a \odot v = (a + v_1, a + v_2, \ldots ,a + v_e)$$
    $$a \odot v \boxplus b \odot w = (\max\{a+v_1,b+w_1\}, \ldots, \max\{a+v_e,b+w_e\}).$$
    \end{definition}

  \begin{example}
  Suppose we have 
  \[
  \begin{array}{rcl}
       v &= & (1, 2, 3),  \\
       w&= & (3, -2, 1), 
  \end{array}
  \]
  and $a = 1, \, b = -3$.  Then we have
   \[
  \begin{array}{rcl}
       a \odot v &= & (1+1, 1+2, 1+3)  \\
       &= & (2, 3, 4), 
  \end{array}
  \] 
  and
   \[
  \begin{array}{rcl}
       a \odot v \boxplus b \odot w &= & (\max\{1+1,(-3)+3\}, \max\{1+2,(-3)+(-2)\}, \max\{1+3,(-3)+1\})  \\
       &= & (2, 3, 4). 
  \end{array}
  \]
  \end{example}

  Throughout this paper we consider the {\em tropical projective torus}, that is, the projective space $\mathbb R^e \!/\mathbb R {\bf
      1}$, where ${\bf 1}:=(1, 1, \ldots , 1)$, the all-one vector.  
  \begin{example}
Consider $\mathbb R^e \!/\mathbb R {\bf 1}$.  Then let
\[
v = (1, 2, 3).
\]
Then over $\mathbb R^e \!/\mathbb R {\bf 1}$ we have the following equality:
\[
v = (1, 2, 3) = (0, 1, 2).
\]
  \end{example}
 
 Note that $\mathbb R^e \!/\mathbb R {\bf 1}$ is isometric to $\mathbb R^{e-1}$.
 
  \begin{example}
Consider $\mathbb R^e \!/\mathbb R {\bf 1}$.  Then let
\[
v = (1, 2, 3), \, w = (1, 1, 1).
\]
Also let $a = -1, \, b = 3$.  Then we have
\[
a \odot v \boxplus b \odot w = (\max(-1 + 1, 3 + 1), \max(-1 + 2, 3 + 1), \max(-1 + 3, 3 + 1)) = (4, 4, 4) = (0, 0, 0).
\]
  \end{example}
  
  In order to conduct a statistical analysis we need a distance measure between two vectors in the space.
  Thus we discuss a distance between two vectors in the tropical projective space.  In fact the following distance between two vectors in the tropical projective space is a metric.  
  \begin{definition}[Generalized Hilbert projective metric]
  For any two points $v, \, w \in \mathbb R^e \!/\mathbb R {\bf 1}$,  the {\em tropical distance} $d_{\rm tr}(v,w)$ between $v$ and $w$ is defined:
    \begin{equation}
  \label{eq:tropmetric} d_{\rm tr}(v,w) \,\, = \,\,
  \max_{i,j} \bigl\{\, |v_i - w_i  - v_j + w_j| \,\,:\,\, 1 \leq i < j \leq e \,\bigr\} = \max_{i} \bigl\{ v_i - w_i \bigr\} - \min_{i} \bigl\{ v_i - w_i \bigr\},
  \end{equation}
  where $v = (v_1, \ldots , v_e)$ and $w= (w_1, \ldots , w_e)$. This
  distance is a metric in $\mathbb R^e \!/\mathbb R {\bf 1}$.  Therefore, we call $d_{\rm tr}$ {\em tropical metric}. 
  \end{definition}
  
  \begin{example}
  Suppose $u_1, \, u_2 \in \mathbb R^3 \!/\mathbb R {\bf 1}$ such that
  \[
  u_1 = (0, 0, 0), \, u_2 = (0, 3, 1).
  \]
  Then the tropical distance between $u_1, \, u_2$ is
  \[
  d_{\rm tr}(u_1, u_2) = \max(0, -3, -1) - \min(0, -3, -1) = 0 - (-3) = 3. 
  \]
  \end{example}
  
  Similar to the BHV metric over the BHV tree space, we need to define a geodesic over the space of ultrametrics.  In order to define a tropical geodesic we need to define a tropical polytope:
  \begin{definition}
  Suppose we have a finite subset $V = \{v_1, \ldots , v_s\}\subset \R^e$ 
  The {\em tropical convex hull} or {\em tropical polytope} of $V$ is the smallest tropically-convex subset containing $V \subset \R^e$ written as the set of all tropical linear combinations of $V$ such that:
    $$
    \mathrm{tconv}(V) = \{a_1 \odot v_1 \boxplus a_2 \odot v_2 \boxplus \cdots \boxplus a_s \odot v_s,
  $$
  where $v_1,\ldots,v_s \in V \mbox{ and } a_1,\ldots,a_s \in \R \}$.
A {\em tropical line segment} between two points $v_1, \, v_2$ is a tropical convex hull of two points $\{v_1, \, v_2\}$.  
  \end{definition}
  Note that the length between two points $u_1, \, u_2 \in \mathbb R^3 \!/\mathbb R {\bf 1}$ along the tropical line segment between $u_1, \, u_2$ equals to the tropical distance $d_{\rm tr}(u_1, u_2)$. 
  In this paper we define a tropical line segment between two points as a {\em tropical geodesic} between these points.  
  
    \begin{example}
  Suppose $u_1, \, u_2 \in \mathbb R^3 \!/\mathbb R {\bf 1}$ such that
  \[
  u_1 = (0, 0, 0), \, u_2 = (0, 3, 1).
  \]
  From the previous example, the tropical distance between $u_1, \, u_2$ is
  \[
  d_{\rm tr}(u_1, u_2) = 3. 
  \]
  Also the tropical line segment between $u_1, \, u_2$ is a line segment between these three points:
  \[
  \begin{array}{c}
        (0, 0, 0)\\
        (0, 2, 0)\\
        (0, 3, 1).\\
  \end{array}
  \]
  The length of the line segment is 
  \[
  d_{\rm tr}((0, 0, 0), (0, 2, 0)) +  d_{\rm tr}((0, 2, 0), (0, 3, 1)) = 2 + 1 = 3.
  \]
  \end{example}
  
  \begin{example}
  Suppose we have a set $ V = \left\{v_1, \, v_2, \, v_3\right\} \subset \mathbb R^3 \!/\mathbb R {\bf 1}$ where
  \[
  v_1 = (0, 0, 0), \, v_2 = (0, 3, 1), \, v_3 = (0, 2, 5).
  \]
  Then we have the tropical convex hull $\mathrm{tconv}(V)$ of $V$ is shown in Figure \ref{fig:tropPoly}.
  \begin{figure}
      \centering
      \includegraphics[width=0.5\textwidth]{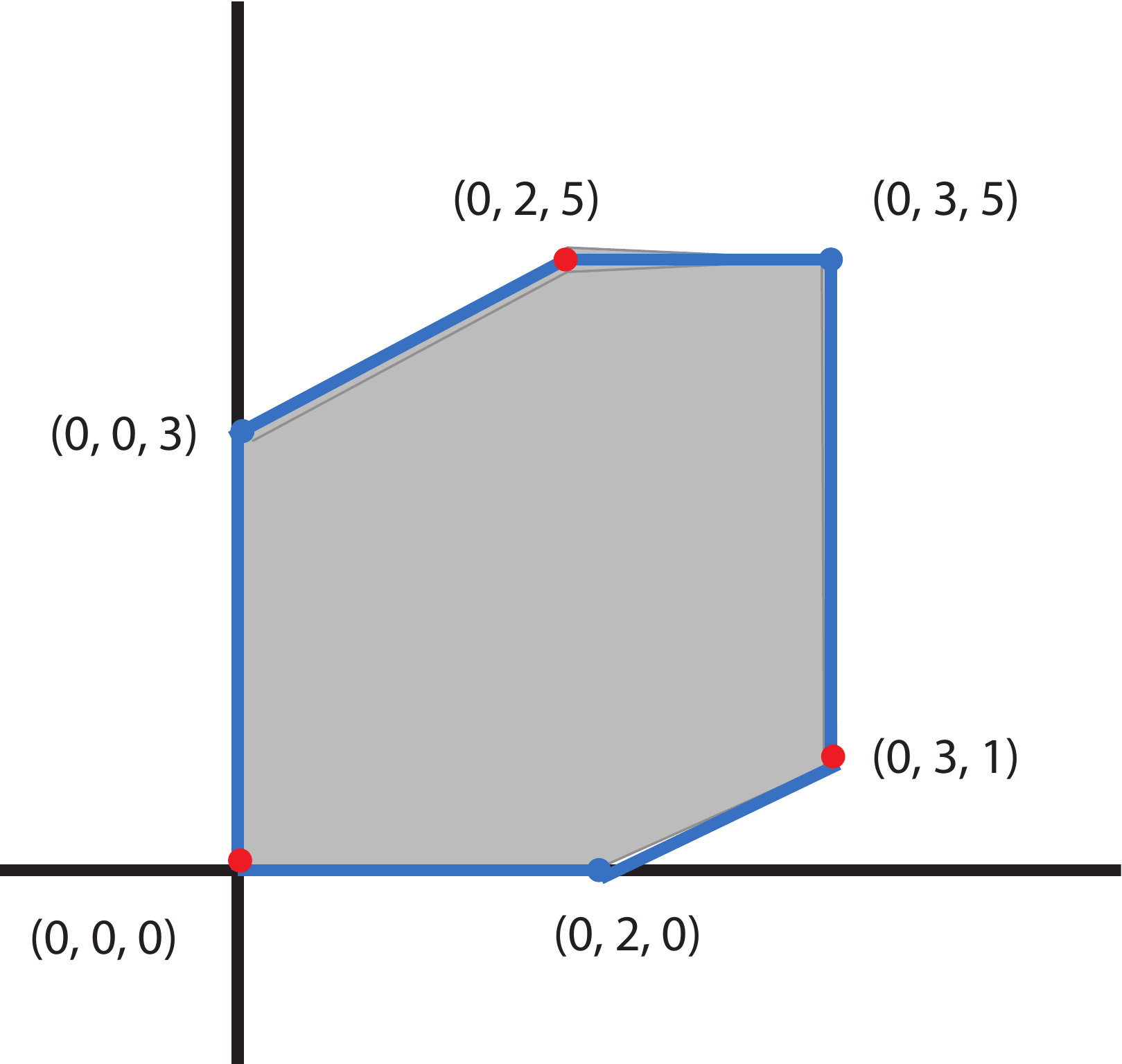}
      \caption{Tropical polytope of three points $(0, 0, 0), \, (0, 3, 1), \, (0, 2, 5)$ in $\mathbb R^3 \!/\mathbb R {\bf 1}$.}
      \label{fig:tropPoly}
  \end{figure}
  \end{example} 

For more details, see the following papers:
\begin{itemize}
    \item D. Maclagan and B. Sturmfels. {\em Introduction to Tropical Geometry} \cite{J}.
\end{itemize}

\section{Tropical Unsupervised Learning}

Unsupervised learning is descriptive and we do not know much about descriptive statistics using tropical geometry with max-plus algebra, for example, tropical Fermat Weber (FW) points and tropical Fr\'ecet means. 

In this section we discuss tropical FW points and tropical Fr\'ecet means, what they are and what we know and we do not know.  In the end of this section, we discuss tropical principal component analysis (PCA).  Over this section we consider the tropical projective torus $\mathbb R^e \!/\mathbb R {\bf 1}$.

\subsection{Tropical Fermat Weber Points}

Suppose we have a sample $\{v_1, \ldots , v_s\}$ over $\mathbb R^e \!/\mathbb R {\bf 1}$.  
A {\em tropical Fermat-Weber point} $y$ minimizes the
sum of distances to the given points. 
\begin{equation}\label{fermat}
y \,:= \, \argmin_{z \in \mathbb R^e \!/\mathbb R {\bf 1}} \sum_{i = 1}^s d_{\rm tr}(z, v_i).
\end{equation}

There are properties of tropical Fermat-Weber points of a sample $\{v_1, \ldots , v_s\}$ over $\mathbb R^e \!/\mathbb R {\bf 1}$.
\begin{proposition}\label{prop:FW1}
Suppose $\mathcal{M} = \RR^e / \RR {\bf 1}$.
Then the set of tropical Fermat-Weber points of a sample $\{v_1, \ldots , v_s\}$ over $\mathbb R^e \!/\mathbb R {\bf 1}$ is a convex polytope.
It consists of all optimal solutions $y = (y_1,\ldots,y_e)$
to the following linear program:
\begin{equation}
\label{eq:ourLP}
 \begin{matrix}
{\rm minimize} \quad d_1 + d_2 + \cdots + d_s \quad \hbox{\rm subject to} \qquad \qquad \qquad \\
\qquad \qquad 
y_j-y_k- v^i_j + v^i_k  \geq -d_i  \quad \hbox{for all} \,\,\, i = 1, \ldots ,s \,\,\hbox{and}\,\,\, 1\le j,k\le e,\\
\qquad \qquad 
y_j-y_k- v^i_j + v^i_k  \leq \phantom{-}d_i 
 \quad \hbox{for all} \,\,\, i = 1, \ldots ,s \,\,\hbox{and} \,\,\, 1\le j,k\le e.
\end{matrix}
\end{equation}
\end{proposition}
From Proposition \ref{prop:FW1}, there can be infinitely many tropical Fermat-Weber points of a sample.  

If we focus on the space of ultrametrics $\mathcal{U}_N$ for {\rm equidistance trees} with $N$ leaves, then we have the following proposition:
\begin{proposition}
If a sample $\{v_1, \ldots , v_s\}$ over the space of ultrametrics $\mathcal{U}_N$, then tropical Fermat-Weber points are in $\mathcal{U}_N$.
\end{proposition}

In \cite{lin2016tropical}, we showed explicitly how to compute the set of all possible Fermat-Weber points in  $\mathbb R^e \!/\mathbb R {\bf 1}$.  However, we do not know  the minimal set of inequalities needed to define the set of all
tropical Fermat-Weber points of a given sample.  Thus here is an open problem:
\begin{problem}
What is the minimal set of inequalities needed to define the set of all
tropical Fermat-Weber points of a given sample?
What is the time complexity to compute the set of tropical Fermat-Weber points
of a sample of $m$ points in $\RR^e/\RR{\bf 1}$?  Is there a
polynomial time algorithm to compute the vertices of the polytope of tropical Fermat-Weber points
of a sample of $s$ points in $\RR^e/\RR{\bf 1}$ in $s$ and $e$?
\end{problem}

For more details, see the following papers:
\begin{itemize}
    \item B. Lin and R. Yoshida {\em Tropical {F}ermat--{W}eber {P}oints} \cite{lin2016tropical}.
\end{itemize}

\subsection{Tropical Fr\'ecet Means}

Suppose we have a sample $\{v_1, \ldots , v_s\}$ over $\mathbb R^e \!/\mathbb R {\bf 1}$.  
A {\em tropical Fr\'echet mean} $y$ minimizes the
sum of distances to the given points. 
\begin{equation}\label{fre}
y \,:= \, \argmin_{z \in \mathbb R^e \!/\mathbb R {\bf 1}} \sum_{i = 1}^s d_{\rm tr}(z, v_i)^2.
\end{equation}

As we formulated computing a tropical Fermat-Weber point as a linear programming problem, we can also formulate computing a tropical Fr\'ecet mean as a quadratic programming problem:
\begin{equation}
\label{eq:ourLP}
 \begin{matrix}
{\rm minimize} \quad d_1^2 + d_2^2 + \cdots + d_s^2 \quad \hbox{\rm subject to} \qquad \qquad \qquad \\
\qquad \qquad 
y_j-y_k- v^i_j + v^i_k  \geq -d_i  \quad \hbox{for all} \,\,\, i = 1, \ldots ,s \,\,\hbox{and}\,\,\, 1\le j,k\le e,\\
\qquad \qquad 
y_j-y_k- v^i_j + v^i_k  \leq \phantom{-}d_i 
 \quad \hbox{for all} \,\,\, i = 1, \ldots ,s \,\,\hbox{and} \,\,\, 1\le j,k\le e.
\end{matrix}
\end{equation}

While we know some propertied of tropical Fermat-Weber points we do not know much about tropical Fr\'echen means.  Here are some basics on tropical Fr\'echet means.
\begin{proposition}
Suppose $\mathcal{M} = \RR^e / \RR {\bf 1}$.
Then the set of tropical Fr\'echen means of a sample $\{v_1, \ldots , v_s\}$ over $\mathbb R^e \!/\mathbb R {\bf 1}$ is a convex polytope.
It consists of all optimal solutions $y = (y_1,\ldots,y_e)$
to the following quadratic program:
\begin{equation}
\label{eq:ourLP}
 \begin{matrix}
{\rm minimize} \quad d_1^2 + d_2^2 + \cdots + d_s^2 \quad \hbox{\rm subject to} \qquad \qquad \qquad \\
\qquad \qquad 
y_j-y_k- v^i_j + v^i_k  \geq -d_i  \quad \hbox{for all} \,\,\, i = 1, \ldots ,s \,\,\hbox{and}\,\,\, 1\le j,k\le e,\\
\qquad \qquad 
y_j-y_k- v^i_j + v^i_k  \leq \phantom{-}d_i 
 \quad \hbox{for all} \,\,\, i = 1, \ldots ,s \,\,\hbox{and} \,\,\, 1\le j,k\le e.
\end{matrix}
\end{equation}
\end{proposition}

Still we do not know much about tropical Fr\'echet means.  First we have the following problem.
\begin{problem}
If a sample $\{v_1, \ldots , v_s\}$ over the space of ultrametrics $\mathcal{U}_N$, then are tropical F\'echet means in $\mathcal{U}_N$?
\end{problem}

We still do not know how to compute tropical Fr\'echet means in efficient ways.  So we have the following problem:
\begin{problem}
Suppose we have $\{v_1, \ldots , v_s\}$ over $\mathbb R^e \!/\mathbb R {\bf 1}$.  Is there an algorithm to compute {\em all} tropical Fr\'echet means in $\mathbb R^e \!/\mathbb R {\bf 1}$?
\end{problem}

\subsection{Tropical Principal Component Analysis (PCA)}

Principal component analysis (PCA) is one of the most popular methods
to reduce dimensionality of input data and to visualize them.  
Classical PCA takes data points in a high-dimensional
Euclidean space and represents them in a lower-dimensional plane in
such a way that the residual sum of squares is minimized.  We cannot directly apply the classical PCA to a set of phylogenetic trees
because the space of phylogenetic trees with a fixed number of leaves is
not Euclidean; it is a union of lower dimensional
polyhedral cones in $\mathbb{R}^{{N \choose 2}}$, where $N$ is the
number of leaves.

There is a statistical method similar to PCA over the space of phylogenetic trees with a fixed set of leaves in terms of the Billera-Holmes-Vogtman (BHV) metric.  

In 2001, Billera-Holmes-Vogtman developed the space of phylogenetic
trees with fixed labeled leaves and they showed that it is ${\rm
  CAT}(0)$ space \cite{BILLERA2001733}.  Therefore, a geodesic between any two points in the space of phylogenetic trees is unique.
  
 Short after that, Nye showed an algorithm in \cite{10.2307/41713594} to compute the first order
principal component 
over the space of phylogenetic trees of $N$ leaves with the BHV metric. 

Nye in \cite{10.2307/41713594} used a convex hull of two points, i.e., the geodesic,
on the tree space as the first order PCA.  However, this idea can not be generalized to higher order principal components with the BHV metric since  the
convex hull of three points with the BHV metric over 
the tree space can have arbitrarily high dimension \cite{doi:10.1137/16M1079841}.  

On the other
hand, the tropical metric in the tree space in terms of the max-plus algebra is well-studied and well-behaved
  \cite{MS}. For example, the dimension of the convex hull of $s$ points in terms of the tropical metric is at most $s-1$. Using the tropical metric, Yoshida et al.\ in \cite{YZZ} introduced a statistical method similar to PCA with the max-plus tropical arithmetic in two ways: the tropical principal linear space, that is, the
best-fit Stiefel
tropical linear space of fixed dimension closest to the data points in
the tropical projective torus; and the tropical principal polytope, that is, the best-fit
tropical polytope with a fixed number of vertices closest to the data
points.  The authors showed that the latter object can be written as a
mixed-integer programming problem to compute them, and they applied the second definition to
datasets consisting of collections of phylogenetic trees. Nevertheless, exactly computing the best-fit tropical polytope can be expensive due to the high-dimensionality of the mixed-integer programming problem.

\begin{definition}
Let $\mathcal P = \tconv(D^{(1)}, \dots, D^{(s)})\subseteq \mathbb
R^{e}/\RR{\bf 1}$ be a tropical polytope with its vertices
$\{D^{(1)},\dots, D^{(s)}\} \subset \mathbb{R}^e/{\mathbb R} {\bf 1}$ and let $S = \{u_1,
\ldots u_n\}$ be a sample from the space of
ultrametrics $\mathcal{U}_N$.
Let $\Pi_{\mathcal P}(S):= \sum_{i = 1}^{|S|}d_{tr}(u_i, u'_i)$, where 
$u'_i$ is the tropical projection of $u_i$ onto a tropical
polytope $\mathcal P$.
Then the vertices $D^{(1)},\dots, D^{(s)}$ of the tropical
polytope $\mathcal P$ are called the {\em $(s-1)$-th order tropical principal 
polytope} of $S$ if the tropical polytope $\mathcal P$ minimizes
$\Pi_{\mathcal P}(S)$ over all possible tropical polytopes with $s$
many vertices.
\end{definition}

In \cite{PYZ}, Page et.al developed a heuristic method to compute tropical principal polytope and they applied it to empirical data sets on genome data of influenza flu collected from New York city, Apicomplexa, and African coelacanth genome data sets.  

Also Page et.al showed the following theorem and lemma:
\begin{theorem}[\cite{PYZ}]
Let $\mathcal P = \tconv(D^{(1)}, \dots, D^{(s)})\subseteq \mathbb R^{e}/\RR{\bf 1}$ be a tropical polytope spanned by ultrametrics in $\mathcal{U}_N$. Then $\mathcal P \subseteq \mathcal{U}_N$ and any two points $x$ and $y$ in the same cell of $\mathcal P$ are also ultrametrics with the same tree topology.
\end{theorem}

\begin{lemma}[\cite{PYZ}]
Let $\mathcal P = \tconv(D^{(1)}, \dots, D^{(s)})\subseteq \mathbb R^{e}/\RR{\bf 1}$ be a tropical polytope spanned by ultrametrics. The origin $\bf 0$ is contained in $\mathcal P$ if and only if the path between each pair of leaves $i,j$ passes through the root of some $D^{(i)}$. 
\end{lemma}

There are still some open problem on tropical PCA.  Here is one of questions we can work on:
\begin{conjecture}
There exists a tropical Fermat-Weber point $x^*\in\mathcal U_N$ of a
sample $D^{(1)}, \dots, D^{(n)}$ of ultrametric trees which is
contained in the $s$th order tropical PCA of the dataset for $s \geq 1$.
\end{conjecture}


For more details, see the following papers:
\begin{itemize}
    \item R. Yoshida, L. Zhang, and X. Zhang. {\em Tropical Principal Component Analysis and its Application to Phylogenetics}. \cite{YZZ}.
    \item R. Page, R. Yoshida, and L. Zhang. {\em Tropical principal component analysis on the space of ultrametrics}. \cite{PYZ}.
\end{itemize}

\section{Tropical Supervised Learning}
For tropical supervised learning, there is not much done.  For classification, there is some work done.  Recently Tang et.al in \cite{TWY} introduced a notion of {\em tropical support vector machines (SVMs)}.  In this section we discuss tropical SVMs and we  introduce a notion of tropical linear discriminant analysis (LDA). 

\subsection{Tropical Classifications}
For tropical classification,  we consider the binary response variables.  Suppose we have a data set given that 
\[
\left\{(x_1, y_1), \ldots , (x_n, y_n)  \right\},
\]
where $x_1, \ldots , x_n \in \mathbb R^{e}/\RR{\bf 1}$ and $y_1, \ldots , y_n \in \{0, 1\}$.  Therefore, the response variable $y_i$ is binary.  Thus, we can partition a sample of data points $x_1, \ldots , x_n \in \mathbb R^{e}/\RR{\bf 1}$ into two sets $P$ and $Q$ such that
\[
\begin{cases}
x_i \in P &\mbox{if } y_i = 0,\\
x_i \in Q & \mbox{if } y_i = 1.\\
\end{cases}
\]

\subsubsection{Tropical support vector machine SVMs}

A support vector machine (SVM) is a supervised learning model to predict the categorical response variable.  For a binary response variable, a classical linear SVM classifies data points by finding a linear hyperplane to separate the data points into two groups.  In this paper we refer a classical SVM as a classical linear SVM over an Euclidean space $\mathbb{R}^e$ with $L_2$ norm.

For an Euclidean space $\mathbb{R}^e$, there are two types of SVMs: {\em hard margin} SVMs and {\em soft margin} SVMs.  A hard margin SVM is a model with the assumption that {\em all} data points can be separated by a linear hyperplane into two groups without errors.  A soft margin SVM is a model which maximizes the margin and also allows some data points in the wrong side of the hyperplane.  

\begin{figure}
      \centering
      \includegraphics[width=0.4\textwidth]{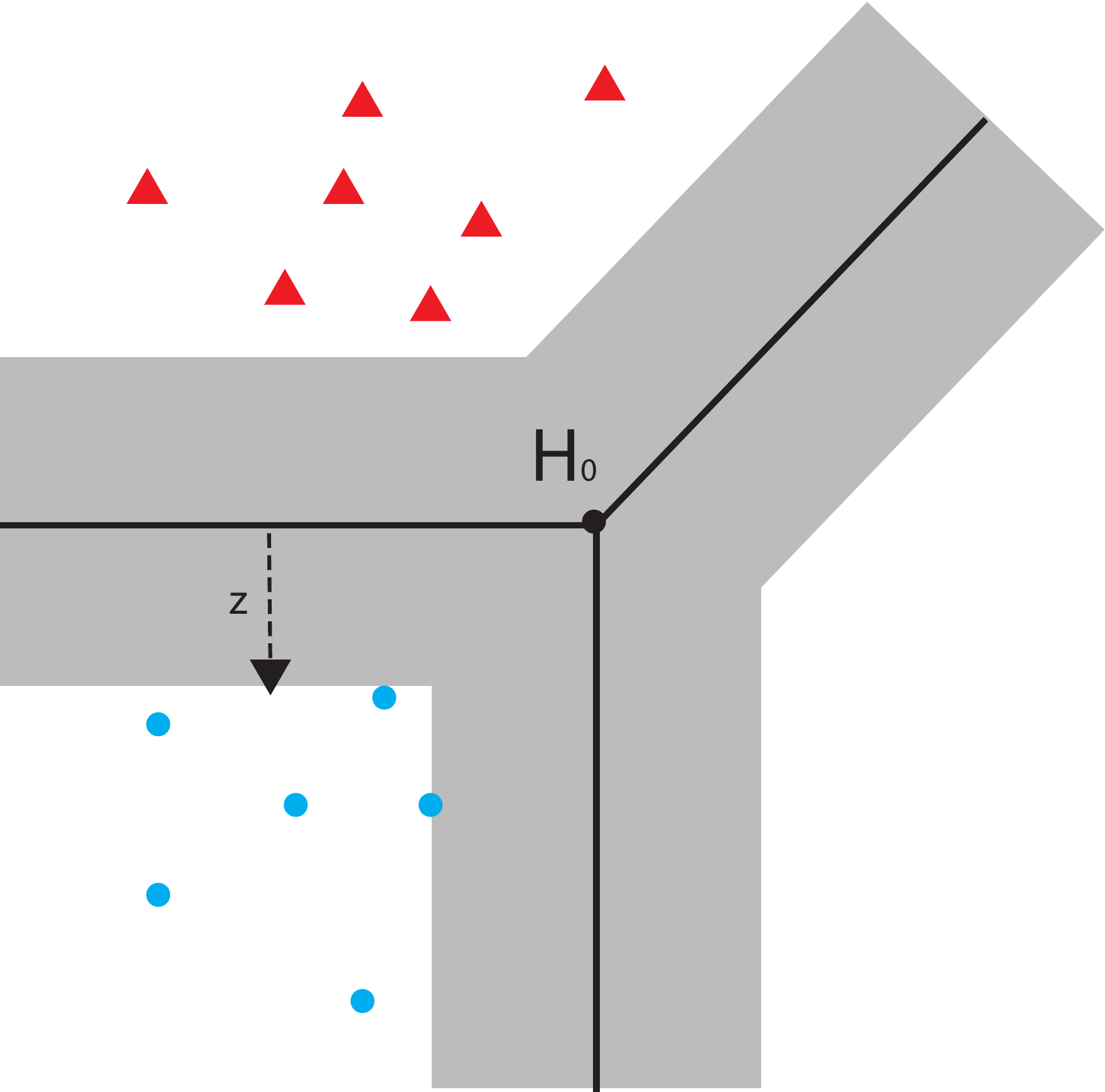}
      \includegraphics[width=0.4\textwidth]{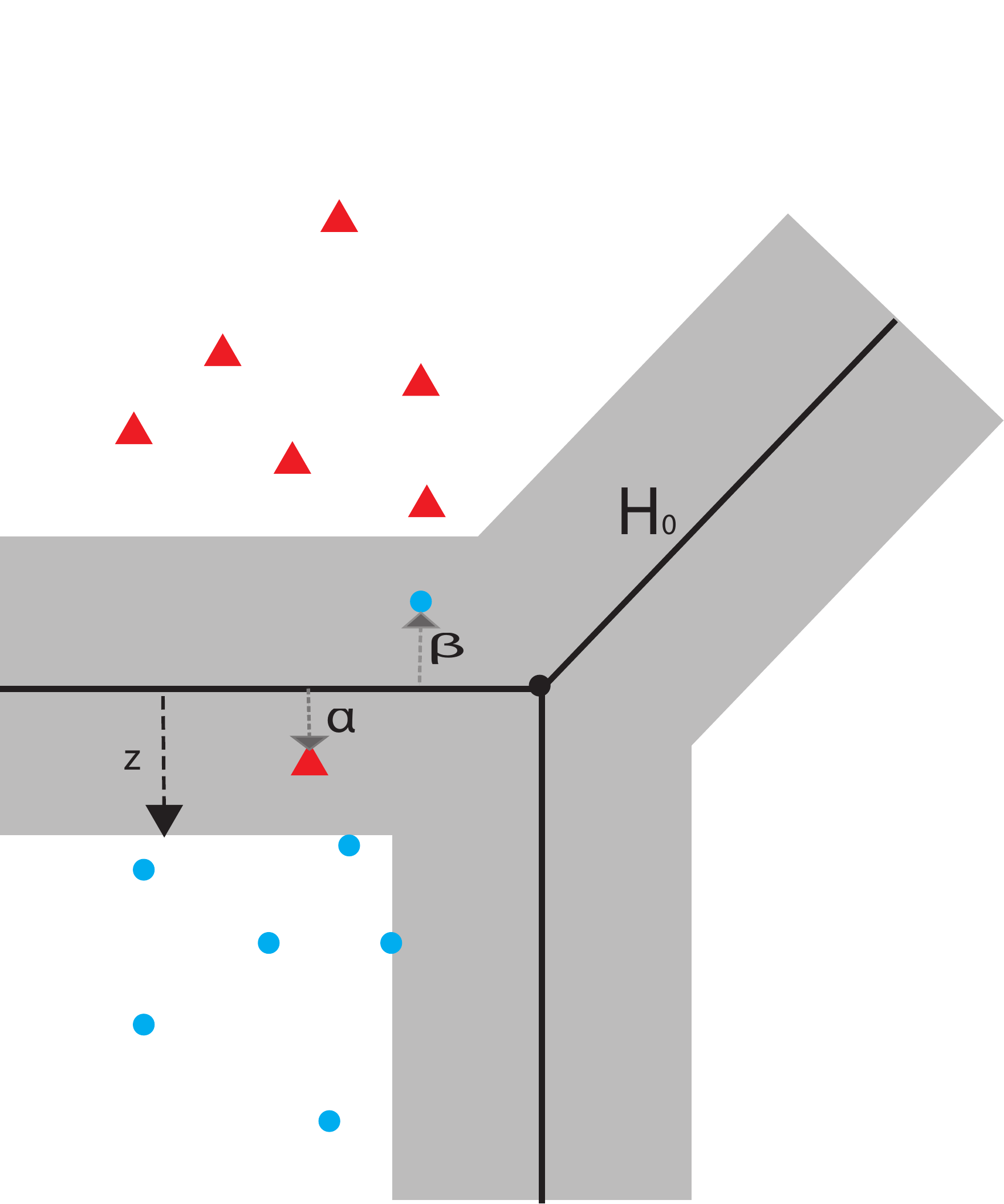}
      \caption{A hard margin tropical SVM (LEFT) and a soft margin tropical SVM (RIGHT) with the binary response variable. A hard margin tropical SVM assumes that all data points from the given sample can be separated by a tropical hyperplane. Red squared dots are data points from $P$ and blue circle dots are data points from $Q$.  A tropical hard margin hyperplane for a tropical hard margin tropical SVM is obtained by maximizing the margin $z$ in the left figure, the distance from the closest data point from the tropical hyperplane (the width of the grey area from the tropical hyperplane in the left figure). A soft margin tropical hyplerplane for a soft margin tropical SVM is obtained by  maximizing a margin similar to a hard margin tropical SVM and by minimizing the sum of $\alpha$ and $\beta$ at the same time.}
      \label{fig:tropHardMargin}
  \end{figure}

Similar to a classical SVM over a Euclidean space, a {\em tropical SVM} is a supervised learning model which classifies data points by finding a tropical hyperplane to separate them.  In \cite{TWY}, as a classical SVM, Tang et.al defined two types of tropical SVMs: {\em hard margin tropical SVMs} and {\em soft margin tropical SVMs}.
A {\em hard margin} tropical SVM introduced by \cite{Gartner} is, similar to a classical hard margin SVM, a model to find a tropical hyperplane which  maximizes the {\em margin}, the minimum tropical distance from data points to the tropical hyperplane (which is $z$ in Figure \ref{fig:tropHardMargin}), to separate these data points into {\em open sectors}.  Note that an open sector of a tropical hyperplane can be seen as a tropical version of an open half space defined by a hyperplane.  A {\em tropical soft margin SVM} introduced by \cite{TWY} is a model to find a tropical SVM to maximizes the margin but it also allows some data points into a wrong open sector.

 The authors in \cite{Gartner} showed that computing a tropical hyperplane for a tropical hard margin SVM from a given sample on the tropical projective space can be formulated as a linear programming problem.  Again, note that, similar to the classical hard margin SVMs, hard margin tropical SVMs assume that there exists a tropical hyperplane such that it separates all data points in the tropical projective space into each open sector (see the left figure in Figure \ref{fig:tropHardMargin}).  

In order to discuss details on tropical SVMs, we need to define a {\em tropical hyplerplane} and their {\em open sectors}.
\begin{definition}
Suppose $\omega:=(\omega_1, \ldots, \omega_e)\in \mathbb R^e \!/\mathbb R {\bf 1}$. The {\em  tropical hyperplane defined by $\omega$}, denoted by $H_{\omega}$, is the set of all points $x\in \mathbb R^e \!/\mathbb R {\bf
    1}$ such that 
  $$\max \{\omega_1+x_1, \ldots , \omega_e+x_e\}$$
    is attained at least twice.  $\omega$ is called the {\em normal vector} of $H_{\omega}$.
\end{definition}

\begin{definition}
A tropical hyperplane $H_{\omega}$ divides the tropical projective space $\mathbb R^e \!/\mathbb R {\bf 1}$ into $e$ components.  These $e$ components divided by $H_{\omega}$ are called {\em open sectors} given that:
  $$S_{\omega}^i~:=~\{\;x\in \mathbb R^e \!/\mathbb R {\bf
    1}\;|\; \omega_i+x_i>\omega_j+x_j,\;\forall j\neq i\;\},\;\;i=1,\ldots,e.$$
\end{definition}

\begin{example}
Consider $\mathbb R^3 \!/\mathbb R {\bf 1}$.  Then a tropical hypoerplane in $\mathbb R^3 \!/\mathbb R {\bf 1}$ has three open sectors seen as Figure \ref{fig:tropHardMargin}.  Note that $\mathbb R^3 \!/\mathbb R {\bf 1}$ is isometric to $\mathbb R^2$. 
\end{example}

Now we define the tropical distance from a point to a tropical hyperplane.
\begin{definition}
  The {\em tropical distance} from a point $x\in \mathbb R^e \!/\mathbb R {\bf
    1}$ to the tropical hyperplane $H_{\omega}$ is defined as:
  $$d_{\rm tr}(x, H_\omega)\;:=\;\min\{d_{\rm tr}(x, y)\;|\;y\in H_{\omega}\}.$$
\end{definition}

A tropical hard margin SVM assumes that all points are separated by a tropical hyperplane and all data points with the same category for their response variable are assigned in the same open sector.  Thus, to compute a tropical hard margin hyperplane for a tropical SVM, we Want to find the normal vector $\omega$ of a tropical  hyperplane $H_{\omega}$  such that 
  \begin{equation*}\label{equation:24}
  \begin{matrix}
  \displaystyle &\max \limits_{\omega \in \mathbb{R}^{d}}\;
  \min \limits_{\xi\in P\cup Q}\;\{\xi_{i(\xi)}+\omega_{i(\xi)}-\xi_{j(\xi)}-\omega_{j(\xi)}\} \\
  \\
  \textrm{s.t.} & \forall \xi \in P\cup Q,\;\forall l \neq i(\xi), j(\xi),\;\; (\xi+\omega)_{l}\leq (\xi+\omega)_{j(\xi)}\leq(\xi+\omega)_{i(\xi)},
  \end{matrix}
  \end{equation*}
  where $i(\xi)$ and $j(\xi)$ are the largest and the second largest coordinate of the vector $\xi+\omega$ for all $\xi \in P\cup Q$.
  
  \begin{theorem}[\cite{TWY}]
  The normal vector $\omega$ of the tropical hard margin $H_{\omega}$ for a tropical SVM is the optimal solution of the following linear programming problem:
      \begin{align}
  &\max \limits_{z \in \mathbb{R}} \; z \label{equation:251} \\
  \textrm{s.t.}\;\;  \forall \xi \in P\cup Q, &\;\;z+\textcolor{black}{\xi_{j(\xi)}}+\omega_{j(\xi)}\textcolor{black}{-\xi_{i(\xi)}}-\omega_{i(\xi)}\leq 0,  \label{equation:252}\\
  \forall \xi \in P\cup Q, &\;\; \omega_{j(\xi)}-\omega_{i(\xi)}\leq \xi_{i(\xi)}-\xi_{j(\xi)}, \label{equation:253} \\
  \forall \xi \in P\cup Q, &\;\forall l \neq i(\xi), j(\xi),\;\; \omega_{l}-\omega_{j(\xi)}\leq \xi_{j(\xi)}-\xi_{l}. \label{equation:254}
  \end{align}
  \end{theorem}

As we discussed earlier, tropical soft margin SVMs are similar to tropical hard margin SVMs.  They try to find a tropical hyperplane which maximizes the margin but also they allow some points to be in a wrong open sector by introducing extra variables $\alpha, \beta$ in Figure \ref{fig:tropHardMargin}.   Tang et.al showed in \cite{TWY} that a soft margin tropical hyperplane for a tropical SVM is the optimal solution of the following linear programming problem such that:
       \begin{align}
  &\max \limits_{\left(z; \alpha; \beta; \gamma \right) \in \mathbb{R}^{(n+2)m+1}} \; z - {\mathcal C}\sum\limits_{\xi\in P\cup Q}\left(\alpha_{\xi}+\beta_{\xi}+\sum\limits_{l\neq i(\xi), j(\xi)}\gamma_{\xi, l}\right)  \label{equation:421} \\
  \textrm{s.t.}\;\;  \forall \xi \in P\cup Q, &\;\;z+\textcolor{black}{\xi_{j(\xi)}}+\omega_{j(\xi)}\textcolor{black}{-\xi_{i(\xi)}}-\omega_{i(\xi)}\leq \alpha_{\xi},  \label{equation:422}\\
  \forall \xi \in P\cup Q, &\;\; \omega_{j(\xi)}-\omega_{i(\xi)}\leq \xi_{i(\xi)}-\xi_{j(\xi)} +\beta_{\xi}, \label{equation:423} \\
  \forall \xi \in P\cup Q, &\;\forall l \neq i(\xi), j(\xi),\;\; \omega_{l}-\omega_{j(\xi)}\leq \xi_{j(\xi)}-\xi_{l} + \gamma_{\xi, l}, \label{equation:424} \\
  \forall \xi \in P\cup Q, &\;\forall l \neq i(\xi), j(\xi),\;\;  \alpha_{\xi}\geq 0, \;\; \beta_{\xi}\geq 0, \;\;\gamma_{\xi, l} \geq 0. \label{equation:425}
  %
\end{align}

There are still many open questions we can ask in terms of tropical SVMs.  In general, if we use methods to find a hard margin or soft margin tropical hyperplane developed in \cite{TWY}, then we have to go through exponentially many linear programming problems. However, we do not know the exact time complexity to find a tropical hard margin or soft margin tropical hyperplane for a tropical SVM.  
\begin{problem}
What is the time complexity of a hard or a soft margin tropical hyperplane for a tropical SVM over the tropical projective torus? Is it NP-hard?
\end{problem}

In addition, the authors in \cite{TWY} focused on tropical hyperplanes for tropical SVMs over the tropical projective torus $\mathbb R^{N \choose 2} \!/\mathbb R {\bf 1}$ not over the space of ultrametrics $\mathcal{U}_N$. 
Again note that $\mathcal{U}_N$ is an union of $N - 1$ dimensional cones over $\mathcal{U}_N \subset \mathbb R^{N \choose 2} \!/\mathbb R {\bf 1}$.  Thus we are interested in how $\mathcal{U}_N$ and a tropical SVM over $\mathbb R^{N \choose 2} \!/\mathbb R {\bf 1}$ related to each other.  More specifically:
\begin{problem}
Can we describe how a hard or soft margin tropical hyperplane for a tropical SVM over the tropical projective torus $\mathbb R^{N \choose 2} \!/\mathbb R {\bf 1}$ separates points in the space of ultrametrics $\mathcal{U}_N$ in terms of geometry?  
\end{problem}

Also  we are interested in defining a tropical SVM over $\mathcal{U}_N$ and developing algorithms to compute them. 
\begin{problem}
Define tropical hard and soft margin "hyperplane" for tropical SVMs over $\mathcal{U}_N$. To define them can we use a tropical polytope instead of a tropical hyperplane?   How can we compute them?  Can we formulate as an optimization problem?
\end{problem}

For more details, see the following papers:
\begin{itemize}
    \item Tang, Wang, and Yoshida.  {\em Tropical Support Vector Machine and its Applications to Phylogenomics}. \cite{TWY}.
\end{itemize}

\subsubsection{Tropical Linear Discriminant Analysis (LDA)}

In this section we discuss {\em tropical linear discriminant analysis} (LDA).  LDA is one of the classical statistical methods to classify dataset into two classes or more as the same time they reduce the dimensionality.  

\begin{figure}
      \centering
      \includegraphics[width=0.9\textwidth]{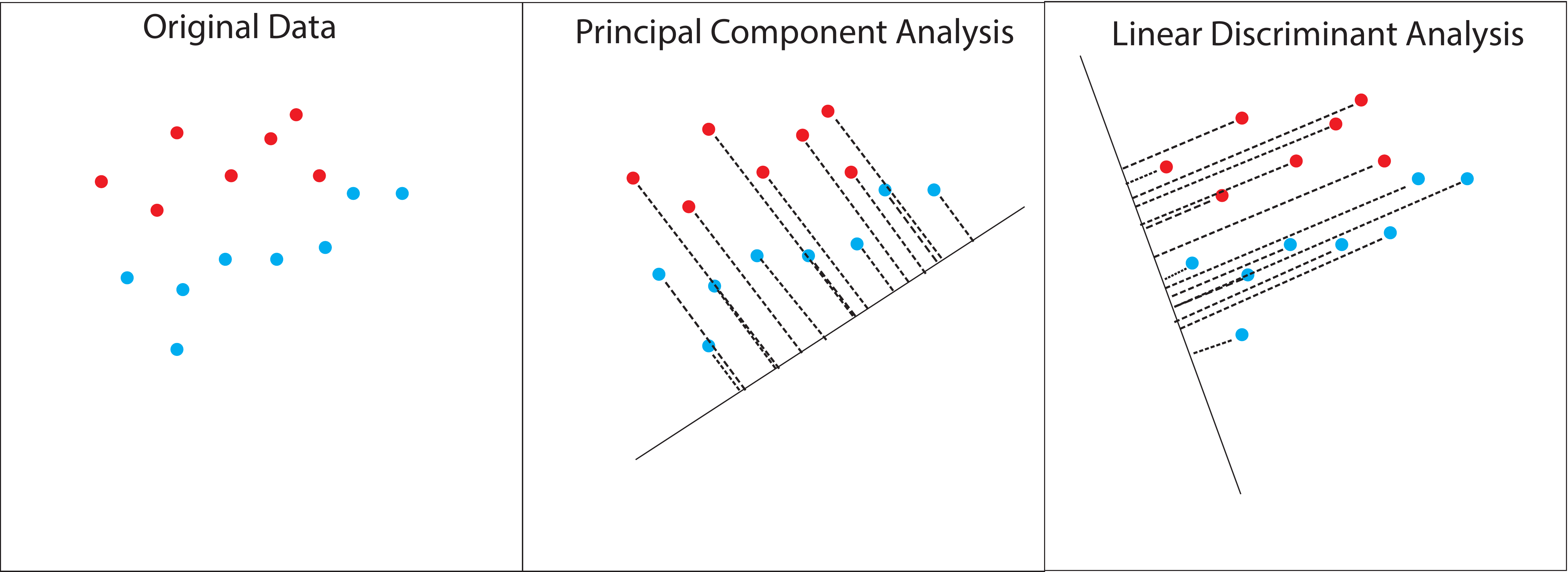}
      \caption{There are two categories in the response variable, red and blue. The middle picture represents PCA and the right picture shows LDA on these points.  }
      \label{fig:LDA}
  \end{figure}

LDA is related to PCA in a Euclidean space and these relations are shown in Figure \ref{fig:LDA}.  The different between PCA and LDA is how to find the direction of a linear plane.  

For two classes of samples $S_1=\{u_1, \ldots , u_{n_1}\}, \, S_2 = \{v_1, \ldots , v_{n_2}\} \subset \mathbb{R}^d$, the linear space for the classical LDA can be found as the optimal solution of an optimization problem such that
\begin{equation}\label{LDA}
    \begin{array}{l}
         \max_w \frac{||\mu_1' - \mu_2'||^2}{{s'}_1^2 + {s'}_2^2} \mbox{ such that } \\
         \\
         \mu_1' =\frac{1}{n_1} \sum_{i = 1}^{n_1} \Proj(u_i),\\
         \\
         \mu_2' =\frac{1}{n_2} \sum_{i = 1}^{n_2} \Proj(v_i),\\
         \\
         {s'}_1 = \sum_{i = 1}^{n_1}||\Proj(u_i) - \mu_1'||,\\
         \\
         {s'}_2 = \sum_{i = 1}^{n_2}||\Proj(v_i) - \mu_2'||,\\
         \\
         \mbox{where }\Proj(\cdot) \mbox{ is a projection onto a linear plane $w$ in } \mathbb{R}^d .
    \end{array}
\end{equation}

Here we use the max-plus algebra in tropical setting.  Also we consider the tropical projective space for now.
Let $d_{\rm tr}$ as a tropical distance between two points in the tropical projective space $\mathbb{R}^d/\mathbb{R} {\bf 1}$.  Then we can formulate the tropical linear space for tropical LDA in Equation \eqref{LDA} as
\begin{equation}\label{troLDA}
    \begin{array}{l}
         \max_w d_{\rm tr}(\mu_1', \mu_2') - {s'}_1 - {s'}_2\mbox{ such that } \\
         \\
         \mu_1' =\arg \min_{z \in w}\sum_{i = 1}^{n_1} d_{\rm tr}(z, \Proj(u_i)),\\
         \\
         \mu_2' =\arg \min_{z \in w}\sum_{i = 1}^{n_2} d_{\rm tr}(z, \Proj(v_i)),\\
         \\
         {s'}_1 = \min_{z \in w}\sum_{i = 1}^{n_1} d_{\rm tr}(z, \Proj(u_i)),\\
         \\
         {s'}_2 = \min_{z \in w}\sum_{i = 1}^{n_2} d_{\rm tr}(z, \Proj(v_i)),\\
         \\
         \mbox{where }\Proj(\cdot) \mbox{ is a projection onto a tropical polytope $w$ in } \mathbb{R}^d/\mathbb{R} {\bf 1} .
    \end{array}
\end{equation}

\begin{problem}
Can we define a tropical LDA over the tropical projective space?  If so how can we find a tropical linear space (or tropical polytope) for a tropical LDA?  
\end{problem}

\begin{problem}
Can we define a tropical LDA over the space of ultrametrics $\mathbb{U}_N$?
\end{problem}

\subsection{Tropical Regression}

For a classical multiple linear regression, with the observed data set
\[
\left\{(x_1, y_1), \ldots (x_n, y_n)\right\}
\]
where $x_i := (x_i^1, \ldots , x_i^e)\in \mathbb{R}^e$ and $y_i \in \mathbb{R}$, 
we try to find a vector $(\beta_0, \beta_1, \ldots , \beta_e) \in \mathbb{R}^{e+1}$ such that
\[
Y = \beta_e X_e + \ldots + \beta_1 X_1 + \beta_0 + \epsilon
\]
where $\epsilon \sim N(0, \sigma)$ with $N(0, \sigma)$ is the Gaussian distribution with the mean $0$ and the standard deviation $\sigma$, $Y$ is a response variable, and $X_1, \ldots , X_e$ are explanatory variables with the smallest following value:
\begin{equation}\label{eq:residual}
 \sum_{i = 1}^n \left(\beta_e x_i^e + \ldots + \beta_1 x_i^1 + \beta_0 - y_i\right)^2 .   
\end{equation}

The value in Equation \ref{eq:residual} is called the sum of squared residuals.  Thus, for a classical multiple linear regression over the Euclidean space $\mathbb{R}^e$, we try to find the linear hyperplane with the smallest sum of squared residuals.  

For tropical regression over the tropical projective space, one can define a tropical regression "polytope" as the tropical polytope with 
\[
\min \sum_{i=1}^n \left(\max \{\beta_e + x_i^e, \ldots , \beta_1 + x_i^1, \beta_0\} - y_i\right)^2 .
\]

It has nothing done in tropical regression.  Thus, it would be interesting to see how one can define them in the tropical projective space as well as the space of ultrametrics.
  \bibliographystyle{plain}
\bibliography{refs}

\end{document}